\documentclass{notemat}

\usepackage{amsmath,amsthm,amssymb}
\usepackage{nmmacro}

\usepackage{url}

\DeclareMathOperator{\tr}{tr}
\DeclareMathOperator{\ad}{ad}
\DeclareMathOperator{\Span}{Span}
\DeclareMathOperator{\Img}{Im}

\def\ip{\langle\cdot,\cdot\rangle}

\begin{document}
\begin{article}
\begin{opening}

\title{On the existence of orthonormal geodesic bases for Lie algebras}

\author{Grant Cairns}
\institute{Department of Mathematics and Statistics, La Trobe University, Melbourne, Australia 3086
\email{G.Cairns@latrobe.edu.au}}

\author{Nguyen Thanh Tung Le\thanks{Part of this work was conducted while the second author was an AMSI summer vacation scholar.}}
\institute{Department of Mathematics and Statistics, La Trobe University, Melbourne, Australia 3086
\email{nt6le@students.latrobe.edu.au}}

\author{Anthony Nielsen}
\institute{Department of Mathematics and Statistics, La Trobe University, Melbourne, Australia 3086
\email{anthony.nielsen@latrobe.edu.au}}

\author{Yuri Nikolayevsky}
\institute{Department of Mathematics and Statistics, La Trobe University, Melbourne, Australia 3086
\email{Y.Nikolayevsky@latrobe.edu.au}}

\runningauthor{Cairns, Le, Nielsen, Nikolayevsky}
\runningtitle{Orthonormal Geodesics Bases}

\begin{abstract}
We show that every unimodular Lie algebra, of dimension at most 4, equipped with an inner product, possesses an orthonormal basis comprised of geodesic elements. On the other hand, we give an example of a solvable unimodular Lie algebra of dimension 5 that has no  orthonormal geodesic basis, for any  inner product. 
\end{abstract}


\end{opening}

\section{Introduction}
 
Let $\mathfrak g $ be a  Lie algebra equipped with an inner product $\ip$. Consider the corresponding simply connected Lie group $G$ equipped with the left invariant Riemannian metric determined by  $\ip$.
A nonzero element $Y\in\mathfrak g $ is said to be a \emph{geodesic vector} if the corresponding left invariant vector field on $G$ is a geodesic vector field. In terms of the Levi-Civita connection $\nabla$, this means that $\nabla_Y Y=0$. This has a simple equivalent expression in terms of the Lie bracket \cite{KNV,KS,Mc,CHGN1}, which we state as a definition.

\begin{Definition} Let $\mathfrak g $ be a  Lie algebra equipped with an inner product $\ip$. Then a nonzero element $Y\in\mathfrak g $ is a \emph{geodesic vector} if $[X,Y]$ is perpendicular to $Y$ for all $X\in\mathfrak g $.
\end{Definition}

Note that some authors use the term  {\em homogeneous geodesic}, to distinguish them from general geodesics on the underlying Lie group. Some authors insist further that $Y$ have unit length. 
For the more general case of totally geodesic subalgebras of Lie algebras, see \cite{Mc,CHGN1,CHGN2}.

\begin{Remark}\label{R} A useful equivalent reformulation of the definition of a geodesic vector is as follows: $Y\in\mathfrak g $ is geodesic if and only if $\Img(\ad(Y))$ is contained in the orthogonal complement of $Y$.
\end{Remark}

Every Lie algebra possesses at least one geodesic vector \cite{Ka,KS,CHGN1}. In \cite{KS} it is shown that semisimple Lie algebras possess an orthonormal basis comprised of geodesics vectors, for every inner product. Results for certain solvable algebras are given in \cite{CKM}. In the present paper we examine the existence of geodesic bases in algebras of low dimension.

\section{Preliminary observations}

First notice that if a Lie algebra $\mathfrak g $, equipped with an orthonormal inner product, has a basis $\{X_1,\ldots,X_n\}$ of geodesic vectors, then for all $i,j$, the element $[X_i,X_j]$ is orthogonal to both $X_i$ and $X_j$. Consequently, the adjoint maps $\ad(X_i)$ have zero trace, and hence $\mathfrak g $ is unimodular. So the natural problem is to determine which unimodular Lie algebras possess an orthonormal basis comprised of geodesic vectors, for some inner product.
We begin with two elementary observations.

\begin{Proposition}\label{P1}
Let $\mathfrak g $ be a nilpotent Lie algebra equipped with an inner product $\ip$. Then there is an orthonormal basis of $(\mathfrak g ,\ip)$ comprised of geodesic vectors.
\end{Proposition}

\begin{pf}
The proof is by induction on the dimension of $\mathfrak g $. The proposition is trivial when $\dim(\mathfrak g )=1$. Suppose that $\dim(\mathfrak g )=n+1$. Since $\mathfrak g $ is nilpotent, its centre is nonzero. Let $Z$ be a central element of length one and set $\mathfrak z  :=\Span(Z)$. Let $\pi: \mathfrak g \to \mathfrak h :=\mathfrak g /\mathfrak z  $ be the natural quotient map and give $\mathfrak h $ the inner product for which the restriction of $\pi$ to the orthogonal complement $\mathfrak z  ^\perp$ of $\mathfrak z  $ is an isometry. By the inductive hypothesis, $\mathfrak h $ possesses an orthonormal basis $\{X_1,\ldots,X_{n}\}$ comprised of geodesic vectors. So, by definition,  for each pair $i,j$, the element $[X_i,X_j]$ is perpendicular to both $X_i$ and $X_j$. For each $i$, let $\bar X_i$ denote the unique element of $\mathfrak z  ^\perp$ with $\pi(\bar X_i)= X_i$. So for each pair $i,j$, the element $[\bar X_i,\bar X_j]$ is perpendicular to both $\bar X_i$ and $\bar X_j$. 
Thus, as $Z$ is central, the elements $\bar X_i$ are geodesic vectors. Furthermore, as $Z$ is central, $Z$ is also geodesic. Hence $\{\bar X_1,\ldots,\bar X_{n},Z\}$
is an orthonormal basis of geodesic  vectors.
\end{pf}

\begin{Proposition}\label{P2}
Let $\mathfrak g $ be a unimodular Lie algebra having a codimension one abelian ideal $\mathfrak h $. Then, for every inner product $\ip$ on $\mathfrak g $,  there is an orthonormal basis of $(\mathfrak g ,\ip)$ comprised of geodesic vectors.
\end{Proposition}

\begin{pf}
We make use of the following  linear algebra result.

\begin{Lemma}[{\cite[Theorem 10]{Pa}}]\label{L}
Suppose that $A$ is a real square matrix with zero trace. Then there is an orthogonal matrix $Q$ such that $QAQ^{-1}$ has zero diagonal.
\end{Lemma}

Let $X\in\mathfrak g $ be a unit vector orthogonal to $\mathfrak h $. Note that $X$ is geodesic since $\mathfrak h $ contains the derived algebra $\mathfrak g '$ of $\mathfrak g $. Choose an orthonormal basis for $\mathfrak h $ and let $A$ denote the matrix representation of the restriction $\ad_{\mathfrak h }(X)$ to $\mathfrak h $ of the adjoint map $\ad(X)$. Since $\mathfrak g $ is unimodular, the matrix $A$ has zero trace. By Lemma \ref{L}, there is an orthonormal basis $\{Y_1,Y_2,\ldots,Y_n\}$  for $\mathfrak h $ relative to which $\ad_{\mathfrak h }(X)$ has zero diagonal. Consequently for each $i$, the element $[X,Y_i]$ is perpendicular to $Y_i$. Thus since $\mathfrak h $ is abelian, the elements $Y_i$ are geodesic  vectors. Hence $\{X,Y_1,Y_2,\ldots,Y_n\}$ is a geodesic basis for $\mathfrak g $.
\end{pf}

\section{Main Results}

It is well known that in dimension less than or equal to three, there are only 5 nonabelian unimodular real Lie algebras; see \cite{TU}. Milnor's classification \cite{Mi} proceeds by considering, for an orthonormal basis $\{X_1,X_2,X_3\}$, the linear map $L(X_i\times X_j):=[X_i,X_j]$, where $\times$ denotes the vector cross product. The matrix of $L$ relative to $\{X_1,X_2,X_3\}$ is symmetric and so its eigenvectors form an orthonormal  basis $\{Y_1,Y_2,Y_3\}$ with, by construction,
\begin{align*}
[Y_2,Y_3]&=\lambda_1 Y_1,\\
[Y_3,Y_1]&=\lambda_2 Y_2,\\
[Y_1,Y_2]&=\lambda_3 Y_3,
\end{align*}
for real coefficients $\lambda_i$. So the basis elements $Y_1,Y_2,Y_3$ are geodesic  vectors.

In dimension 4, there are infinitely many isomorphism classes of unimodular  real Lie algebras; see \cite{Mu,deG,St}. Instead of using the classification, we will argue directly. The following result resolves a question raised in  \cite{Ni}.

\begin{Theorem}\label{T}
Let $\mathfrak g $ be a unimodular Lie algebra of dimension 4 equipped with an inner product $\ip$. Then there is an orthonormal basis of $(\mathfrak g ,\ip)$ comprised of geodesic  vectors.
\end{Theorem}

\begin{pf}
Let $\mathfrak g $ be as in the statement of the theorem. First observe that if $\mathfrak g $ is not solvable, then from Levi's Theorem and the classification of semisimple algebras \cite{Ja}, there is a Lie algebra isomorphism $\mathfrak g \cong \mathfrak h \oplus\mathbb R $ where $\mathfrak h $ is a simple Lie algebra of dimension three; so $\mathfrak h $ is isomorphic to either $\mathfrak{so}(3,\mathbb R )$ or $\mathfrak{sl}(2,\mathbb R )$ though we won't need this fact. Let $W\in\mathfrak g $ be a unit vector orthogonal to $\mathfrak h $. Since $\mathfrak h =\mathfrak g '$, the element $W$ is geodesic. Note that $W$ can be written uniquely in the form $W=X+Z$, where $Z$ is in the centre of $\mathfrak g $ and $X\in\mathfrak h $. As we discussed above,  $\mathfrak h $ has an orthonormal basis $\{Y_1,Y_2,Y_3\}$  of elements that are geodesic  vectors of $\mathfrak h $. Then for each $i=1,2,3$,
\[
[W,Y_i]=[X,Y_i]\in \Img(\ad(Y_i)|_{\mathfrak h }),
\]
which is perpendicular to $Y_i$ by Remark \ref{R}. Hence $Y_i$ is geodesic in $\mathfrak g $ and thus $\{Y_1,Y_2,Y_3,W\}$ is an orthonormal basis of geodesic vectors. So we may assume that $\mathfrak g $ is solvable.

Consider the derived algebra $\mathfrak g ':=[\mathfrak g ,\mathfrak g ]$. As $\mathfrak g $ is solvable, $\mathfrak g '$ is nilpotent \cite[Chap.~1.5.3]{Bou}. If $\dim(\mathfrak g ')=0$, the algebra $\mathfrak g $ is abelian and the theorem holds trivially. We consider three cases according to the remaining possibilities for the dimension of $\mathfrak g '$.

If $\dim(\mathfrak g ')=1$, let $\mathfrak g '=\Span(W)$. For each $X\in \mathfrak g $ we have $[X,Y]\in \mathfrak g '$ for all $Y\in\mathfrak g $ and hence $[X,W]=\tr(\ad(X))W=0$. Thus $\mathfrak g '$ is contained in the centre of $\mathfrak g $. Consequently $\mathfrak g $ is nilpotent, and the theorem holds by Proposition \ref{P1}.

If $\dim(\mathfrak g ')=3$, then either $\mathfrak g '$ is abelian or $\mathfrak g '$ is isomorphic to the Heisenberg Lie algebra $\mathfrak h _1$. If $\mathfrak g '$ is abelian, the required result follows from Proposition \ref{P2}. So suppose that $\mathfrak g '$ is isomorphic to the Heisenberg Lie algebra. Choose an orthonormal basis $\{X,Y,Z\}$ of $\mathfrak g '$ with $Z$ in the centre of $\mathfrak g '$ and $[X,Y]=\lambda Z$ for some $\lambda\not=0$. Let $W\in\mathfrak g $ be a unit vector orthogonal to $\mathfrak g '$. Note that $W$ is geodesic. Let $A$ denote the matrix representation of the restriction $\ad_{\mathfrak g '}(W)$ to $\mathfrak g '$ of the adjoint map $\ad(W)$. Since the centre of $\mathfrak g '$ is a characteristic ideal \cite[Chap.~1.1.3]{Bou}, it is left invariant by $\ad(W)$. Hence $A$ has the form
\[
A=\begin{pmatrix}
a&b&0\\
c&d&0\\
e&f&g\end{pmatrix},
\]
for reals $a,b,c,d,e,f,g$.  By the Jacobi identity,
\[
[\ad(W)X,Y]+[X,\ad(W)Y]=\ad(W)\lambda Z=g\lambda Z.
\]
Thus
\[
[aX+cY+eZ,Y]+[X,bX+dY+fZ]=(a+d)\lambda Z=g\lambda Z,\]
and so $a+d=g$. Since $\ad(W)$ has zero trace, $a+d+g=0$, and hence $a+d=g=0$. So, by Lemma \ref{L}, we can make an orthonormal change of basis for the space $\Span(X,Y)$ so that relative to the new basis, $a=d=0$. Then  $\{W,X,Y,Z\}$ is an orthonormal geodesic basis for $\mathfrak g $.

Finally, if $\dim(\mathfrak g ')=2$, then $\mathfrak g '$ is abelian. Let  $\{X,Y\}$ be an orthonormal basis for the orthogonal complement $\mathfrak g '^\perp$ to $\mathfrak g '$. 
Note that $[X,Y]\in\mathfrak g '$ and so as $\mathfrak g '$ is abelian, the adjoint map $\ad([X,Y])$ acts trivially on $\mathfrak g '$. Thus, since $\ad(X)\circ \ad(Y)-\ad(Y)\circ \ad(X)=\ad([X,Y])$ by the Jacobi identity, we have that the restrictions $\ad_{\mathfrak g '}(X)$ and $\ad_{\mathfrak g '}(Y)$ commute. Since $\ad_{\mathfrak g '}(X)$ and $\ad_{\mathfrak g '}(Y)$ have zero trace,  relative to any choice of orthonormal basis for $\mathfrak g '$, they have matrix representations in $\mathfrak{sl} (2,\mathbb R )$. However, it is well known and easy to see that the maximum abelian subalgebras of  $\mathfrak{sl} (2,\mathbb R )$ have dimension one.
Consequently, by orthonormal change of basis for $\mathfrak g '^\perp$ if necessary, we may assume that $\ad_{\mathfrak g '}(Y)\equiv 0$. Then $\Span(Y,\mathfrak g ')$ is a codimension one abelian ideal and the result follows from Proposition \ref{P2}.
\end{pf}

Note that in Propositions \ref{P1}, \ref{P2} and Theorem \ref{T}, the required orthonormal geodesic basis is obtained for \emph{every}  inner product. The following example gives a solvable unimodular Lie algebra that has no  orthonormal geodesic basis for \emph{any}  inner product.

\begin{Example} \label{E}
Consider the 5-dimensional algebra $\mathfrak g $ with basis $B=\{X_1,\ldots,X_5\}$ and relations
\[
  \begin{array}{lll}
  & [X_1,X_2] = 3X_2  & [X_2,X_3] = X_4 \\
  & [X_1,X_3] = -4X_3 & [X_2,X_4] = X_5 \\
  & [X_1,X_4] = -X_4  & \\
  & [X_1,X_5] = 2X_5. &
  \end{array}
\]
The ideal $\mathfrak n$ generated by $X_2,\ldots,X_5$ is the (unique) 4 dimensional filiform nilpotent Lie algebra. The adjoint map $\ad(X_1)$ is obviously a Lie derivation of $\mathfrak n$, so the Jacobi identities of $\mathfrak g $ hold. Clearly, $\mathfrak g $ is solvable and unimodular, but not nilpotent. We will show that $\mathfrak g $ has no orthonormal basis  comprised of geodesic vectors, for any inner product. Equip $\mathfrak g $ with an inner product and suppose it has a basis (not necessarily orthonormal) of geodesic elements $Y_1,\ldots,Y_5$.
Consider an arbitrary geodesic   vector $Y = a_1X_1+\cdots+a_5X_5 \in \mathfrak g $. Let $V_i:=\Span(X_i,X_{i+1},\ldots,X_5)$. Relative to the basis $B$, the matrix representation of $\ad(Y)$ is
\begin{equation} \label{GV_eq3}
A=\begin{pmatrix}
0     & 0    & 0      & 0      & 0      \\
-3a_2 & 3a_1 & 0 & 0 & 0 \\
4a_3  & 0    & -4a_1  &   0     & 0       \\
a_4   & -a_3 & a_2    & -a_1   &   0     \\
-2a_5 & -a_4 & 0      & a_2    & 2a_1
\end{pmatrix}
.
\end{equation}
A priori, there are five cases:
\begin{enumerate}
\item\label{i1} $Y\in V_5$; that is, $a_1 = a_2 = 0,a_3 = 0,a_4 = 0,a_5 \neq 0$,
\item\label{i2} $Y\in V_4\backslash V_5$; that is, $a_1 = a_2 = 0,a_3 = 0,a_4 \neq 0$,
\item\label{i3} $Y\in V_3\backslash V_4$; that is, $a_1 = a_2 = 0,a_3 \neq 0$,
\item\label{i4} $Y\in V_2\backslash V_3$; that is, $a_1 = 0,a_2 \neq 0$,
\item\label{i5} $Y\in V_1\backslash V_2$; that is, $a_1 \neq 0$.
\end{enumerate}
In fact, case \eqref{i1} is impossible as otherwise $Y$ is a nonzero multiple of $X_5$ but $\Img(\ad(Y))=V_5$, which must be orthogonal to $Y$ by Remark \ref{R}.
Similarly, case \eqref{i2} is impossible as otherwise $Y\in V_4$ but $\Img(\ad(Y))=V_4$, which has to be orthogonal to $Y$ by Remark \ref{R}.

In case \eqref{i3}, $\Img(\ad(Y))=\Span(X_3+\frac{a_4}{4a_3}X_4-\frac{2a_5}{4a_3}X_5,X_4+\frac{a_4}{a_3}X_5)$.

In case \eqref{i4}, $\Img(\ad(Y))=\Span(X_2-\frac{4a_3}{3a_2}X_3,X_4,X_5)$ and in particular, $Y$ is orthogonal to $V_4$.

In case \eqref{i5},  $\Img(\ad(Y))=V_2$, which is the orthogonal complement of $Y$ by Remark \ref{R}. It follows that, up to a constant multiple, there is at most one geodesic vector $Y$ with $a_1\not=0$.  Thus, there is at most one basis element,  $Y_1$ say, with $Y_1\not\in V_2$. 

So we have established that 
\[
Y_2,Y_3,Y_4,Y_5\in  V_2\backslash V_4.
\]
But if $Y_i\in  V_2\backslash V_3$, then as we saw in case~\eqref{i4}, $Y_i$  is orthogonal to $V_4$. So for dimension reasons, there are at most two of the $Y_i$  in $V_2\backslash V_3$. 

Now suppose the $Y_i$ are orthogonal. If two of the $Y_i$ are in $V_2\backslash V_3$,  they would both be orthogonal to $V_4$. This would force the remaining two $Y_i$ to be in $V_4$, which we have seen is impossible. So only one  of the $Y_i$, say $Y_2$, can be  in $V_2\backslash V_3$. So $Y_3,Y_4,Y_5$ are in $V_3\backslash V_4$. Note that $V_3$ is left invariant by $\ad(X_1)$. Relative to the orthonormal basis $\{Y_3,Y_4,Y_5\}$ for $V_3$, the map $f:=\ad(X_1)|_{V_3}$ has zero diagonal, because the elements $Y_3,Y_4,Y_5$ are geodesic vectors. So $f$ has zero trace. But relative to the basis $X_3,X_4,X_5$, it is clear that $f$ has trace $-4-1+2=-3$. This is a contradiction.
\end{Example}

\begin{Remark}
Although the algebra of the above example does not posses an orthonormal basis of geodesic vectors, it does posses an inner product for which there is a (nonorthonormal) basis of geodesic vectors. That is, we claim that there exists an inner product such that the span of all the geodesic vectors is the whole algebra. Using the above notation, regardless of the choice of an inner product, there is exactly one geodesic vector $Y_1$ (up to scaling) not lying in $\mathfrak g '=V_2$, namely any nonzero vector from $(\mathfrak g ')^\perp$. We therefore want to show that the span of all the geodesic vectors from $\mathfrak g '$ covers $\mathfrak g '$. We will now specify the inner product. First choose $X_2, X_3$ to be orthonormal and orthogonal to $V_4$. Then from case~\eqref{i4}, a vector $Y=a_2X_2+a_3X_3+a_4X_4+a_5X_5$ (with $a_2 \ne 0$) is geodesic if it is orthogonal to $V_4$ (so $a_4=a_5=0$) and $\sqrt{3}a_2=\pm2a_3$. This gives two linearly independent geodesic vectors $Y_2, Y_3 \in \mathfrak g ' \cap V_4^\perp$, neither of which lies in $V_3$. It remains to show that there exist at least two geodesic vectors in $V_3$ whose projections to $V_4$ are linearly independent. Define the remaining components of the inner product by requiring that $\|X_4\|=\|X_5\|=1$ and $\langle X_4,X_5\rangle=\varepsilon \in (0,1)$. Then by case~\eqref{i3}, a vector $Y=a_3X_3+a_4X_4+a_5X_5$ (with $a_3 \ne 0$) is geodesic if and only if the following two conditions hold.
\begin{align} \label{eq:1}
    4a_3^2+a_4^2-2a_5^2-\varepsilon  a_4a_5 &=0, \\
    a_3(a_4+\varepsilon  a_5) +a_4(\varepsilon  a_4+a_5)&=0. \label{eq:2}
\end{align}
Note that $a_4+\varepsilon  a_5 \ne 0$, as otherwise from~\eqref{eq:2} either $a_4=0$ or $\varepsilon   a_4+a_5=0$ and in both cases we obtain $a_4=a_5=0$ and then  \eqref{eq:1} would give $a_3=0$, which is a contradiction. Solving~\eqref{eq:2} for $a_3$ and substituting to~\eqref{eq:1} we get
\begin{equation}\label{eq:3}
    4a_4^2(\varepsilon  a_4+a_5)^2+(a_4^2-2a_5^2-\varepsilon  a_4a_5)(a_4+\varepsilon  a_5)^2=0.
\end{equation}
Note that as $a_3\not=0$ and $a_4+\varepsilon  a_5 \ne 0$, (\ref{eq:2}) gives $a_4 \ne 0$. Dividing \eqref{eq:3}
by $a_4^4$ and taking $t=a_5a_4^{-1}$ we obtain 
\[
-2\varepsilon ^2t^4-(4\varepsilon +\varepsilon ^3)t^3+(2-\varepsilon ^2)t^2+9\varepsilon  t+(1+4\varepsilon ^2)=0.
\]
The polynomial on the left-hand side has at least one positive root $t_+$ and at least one negative root $t_-$. It follows that $V_3$ contains at least two geodesic vectors $Y_{\pm}=(a_3)_{\pm} X_3 + X_4 + t_{\pm} X_5$ (where $(a_3)_{\pm}$ are determined by~\eqref{eq:2}). Their projections to $V_4$ are linearly independent, hence $\Span(Y_1,Y_2,Y_3,Y_+,Y_-)=\mathfrak g $.

Note that for an arbitrary inner product on this algebra, the geodesic vectors may not span the entire algebra. For example, choosing an inner product with all the $X_i$'s orthogonal we obtain that the span of all the geodesic vectors is the proper subspace $X_4^\perp$ of $\mathfrak g $.
\end{Remark}

We conclude this paper with two questions that have arisen from this study:

\medskip
\noindent
{\bf Question 1}. Is it true that every unimodular Lie algebra possesses an inner product for which the geodesic vectors span the algebra?

\medskip
\noindent
{\bf Question 2}. Apart from nilpotent Lie algebras, are there natural families of unimodular Lie algebras that possess an inner product for which there is an orthonormal basis of  geodesic vectors?

\begin{acknowledgements}
The authors are very grateful to Ana Hini\'c Gali\'c who detected a number of typos in the first version of this paper.
\end{acknowledgements}


\providecommand{\bysame}{\leavevmode\hbox to3em{\hrulefill}\thinspace}
\providecommand{\MR}{\relax\ifhmode\unskip\space\fi MR }
\providecommand{\MRhref}[2]{%
  \href{http://www.ams.org/mathscinet-getitem?mr=#1}{#2}
}
\providecommand{\href}[2]{#2}

\end{article}
\end{document}